\documentclass[12pt,reqno]{amsart}
\usepackage{amssymb}
\usepackage{mathrsfs}
\usepackage{a4}
\usepackage[all]{xy}

\newcommand{\heute}{24 January 2009}

\theoremstyle{plain}
\newtheorem{theorem}{Theorem}[section]

\newtheorem{lemma}[theorem]{Lemma}
\newtheorem{corollary}[theorem]{Corollary}

\theoremstyle{remark}

\newtheorem*{defn}{Definition}
\newtheorem*{rk}{Remark}

\newtheorem*{notn}{Notation}

\newcommand{\dashTwo}[1]{\textup{(\ref{two}${}'$)}}

\newcommand{\ignore}[1]{}

\newcommand{\f}[1][p]{\mathbb{F}_{#1}}

\newcommand{\Gro}[1]{Gr\"ob\-ner}

\newcommand{\E}{\mathscr{E}}
\newcommand{\pp}[1][G,V]{\mathscr{P}(#1)}

\newcommand{\abs}[1]{\left|#1\right|}
\newcommand{\scr}[1]{\mathscr{#1}}

\begin{document}

\title{On Oliver's $p$-group conjecture: II}
\author[D.~J. Green]{David J. Green}
\address{Department of Mathematics \\
Friedrich-Schiller-Universit\"at Jena \\ 07737 Jena \\ Germany}
\email{David.Green@uni-jena.de}
\author[L.~H\'ethelyi]{L\'aszl\'o H\'ethelyi}
\address{Department of Algebra \\
Budapest University of Technology and Economics \\ Budapest \\ Hungary}
\thanks{H\'ethelyi supported by the Hungarian Scientific Research Fund
OTKA, grant T049 841}
\author[N.~Mazza]{Nadia Mazza}
\address{Department of Mathematics and Statistics \\ Fylde College \\
Lancaster University \\ Lancaster LA1 4YF \\ United Kingdom}
\date{\heute}

\begin{abstract}
Let $p$ be an odd prime and $S$ a finite $p$-group. B.~Oliver's conjecture
arises from an open problem in the theory of $p$-local finite groups. It is
the claim that a certain characteristic subgroup
$\mathfrak{X}(S)$ of $S$ always contains the Thompson subgroup.
In previous work the first two authors
and M.~Lilienthal recast Oliver's conjecture
as a statement about the representation theory of the factor
group $S/\mathfrak{X}(S)$. We now verify the conjecture for a wide variety
of groups~$S/\mathfrak{X}(S)$.
\end{abstract}

\maketitle

\section{Introduction}
\label{intro}
\noindent
Let $p$ be an odd prime and $S$ a finite $p$-group. 
An open question in the theory of $p$-local finite groups asks
whether there is a unique centric linking system associated to
every fusion system
(see the survey article~\cite{BLO:survey} by Broto, Levi and Oliver)\@.
Bob Oliver derived in~\cite{Oliver:MartinoPriddyOdd} a purely group-theoretic
conjecture which would imply existence and uniqueness of the linking system,
at least at odd primes.
He constructed a characteristic subgroup $\mathfrak{X}(S)$,
and conjectured that it
always contains the Thompson subgroup $J(S)$ generated
by the elementary abelian subgroups of greatest rank. The first two authors
and M.~Lilienthal studied Oliver's conjecture in~\cite{oliver}
and recast it as a question
about the quotient group $G = S / \mathfrak{X}(S)$.

In this paper we shall use
methods from the area of finite group theory known as Thompson factorization
(see \S32 in Aschbacher's book~\cite{Aschbacher:book}) to study the
properties of certain faithful $\f G$-modules which arise in this reformulation of
Oliver's conjecture. This allows us to prove the conjecture for a wide
variety of quotient groups~$G$. Our main result is as follows:

\begin{theorem}
\label{thm:main}
Suppose that $p$~is an odd prime and $S$ is a $p$-group such that
$S / \mathfrak{X}(S)$ satisfies any of the following conditions
\begin{enumerate}
\item
its (nilpotence) class is at most four;
\item
it is metabelian;
\item
it is of maximal class;
\item
its $p$-rank is at most~$p$.
\end{enumerate}
Then Oliver's conjecture $J(S) \leq \mathfrak{X}(S)$ holds for~$S$.
\end{theorem}

\begin{rk}
Alperin showed that every regular $3$-group is
metabelian~\cite{Alperin:Regular}\@. So Oliver's conjecture also
holds if $S / \mathfrak{X}(S)$ is a regular $3$-group.
\end{rk}

\begin{rk}
Every $p$-group $G$ does indeed occur as such a quotient
$G = S/\mathfrak{X}(S)$, by Lemma~2.3 of~\cite{oliver}\@.
A computation using GAP~\cite{GAP4} shows that the iterated wreath
product group $G = C_3 \wr C_3 \wr C_3$ satisfies none of the above conditions.
This group has order $3^{13}$ and class~$9$, so it is neither class${}\leq4$
or maximal class. The derived subgroup has class~$3$, so it is not metabelian.
And the ``double diagonal'' subgroup shows that the rank is at least~$9$.
\end{rk}

\noindent
Theorem~\ref{thm:main} follows by translating the following two module-theoretic
results back into the original language of Oliver's conjecture.
The notions of ``quadratic element'' and ``$F$-module'' are recalled in
\S\ref{section:notation} and \S\ref{section:offenders}, respectively.

\begin{theorem}
\label{thm:main1}
Suppose that $G$ is a $p$-group and $V$ is a faithful $\f G$-module such that
$\Omega_1(Z(G))$ has no quadratic elements. If $G$ satisfies either of the
following conditions
\begin{enumerate}
\item
$G$ has class at most four;
\item
$G$ is metabelian;
\end{enumerate}
then $V$~cannot be an $F$-module.
\end{theorem}

\noindent
Notice that the assumption that there are no quadratic elements  in
$\Omega_1(Z(G))$ implies that the prime~$p$ has to be odd. Recall now that
the rank of an elementary abelian $p$-group is its dimension as
$\f$-vector space, and that the $p$-rank of a finite group is the
maximum of the ranks of its elementary abelian $p$-subgroups.

\begin{theorem}
\label{thm:main2}
Suppose that $p$~is an odd prime, that $G$ is a $p$-group, and that $V$ is
a faithful $\f G$-module such that every non-identity element of $\Omega_1(Z(G))$
acts with minimal polynomial $X^p-1$. If $G$ satisfies either of the
following conditions
\begin{enumerate}
\item
$G$ is of maximal class;
\item
the $p$-rank of $G$ is at most~$p$;
\end{enumerate}
then $V$~cannot be an $F$-module.
\end{theorem}

\noindent
A key step in the proof of Theorem~\ref{thm:main1} is the following result.
We recall the term ``offender'' in \S\ref{section:offenders} below.

\begin{theorem}
\label{thm:normalAbelian}
Let $G$ be a $p$-group and $V$ a faithful $\f G$-module such that there
are no quadratic elements in $\Omega_1(Z(G))$.
\begin{enumerate}
\item
\label{enum:Enormal}
If $A$ is an abelian normal subgroup of $G$, then $A$ does not contain
any offender.
\item
\label{enum:Gdash}
Suppose that $E$~is an offender. Then $[G',E] \neq 1$.
\end{enumerate}
\end{theorem}

\begin{proof}[Proof of Theorem~\ref{thm:main}]
The condition in Theorem~\ref{thm:main2} on the minimal polynomial of
each non-identity element of $\Omega_1(Z(G))$ is condition~(PS)
of~\cite{oliver}\@. 
Since $p > 2$, this condition implies that there
are no quadratic elements in $\Omega_1(Z(G))$. Now combine our
Theorems \ref{thm:main1}~and \ref{thm:main2} with
Theorem~1.2 of~\cite{oliver}\@.
\end{proof}

\begin{proof}[Proof of Theorem~\ref{thm:main1}]
Theorem~\ref{thm:class4} deals with the case of class
at most four, and Theorem~\ref{thm:abelianDerived} treats
the metabelian case.
\end{proof}

\begin{proof}[Proof of Theorem~\ref{thm:main2}]
Theorem~\ref{theorem:rank-p} handles the case $p$-rank$(G) \leq p$.
Lemma~\ref{lemma:maxlClassRankp} shows that groups of maximal
class have $p$-rank at most~$p$.
\end{proof}

\noindent
For sake of completeness, we remind the reader that the Oliver
subgroup $\mathfrak{X}(S)$ is the
largest normal subgroup of $S$ that has a $Q$-series; that is, there are
an integer $n\geq1$ and a series
$1=Q_0\leq Q_1\leq\dots\leq Q_n=\mathfrak{X}(S)$
of normal subgroups $Q_i$ of $S$ and such that
$[\Omega_1(C_S(Q_{i-1})),Q_i;p-1]=1$, for all $1\leq i\leq n$.

\emph{Structure of the paper} \quad
In \S\ref{section:notation} we comment on a few commutator relations in
a semidirect product and
introduce some handy notation for these. In
the short~\S\ref{section:MeierStell} we recall a lemma of Meierfrankenfeld
and Stellmacher. We prove Theorem~\ref{thm:normalAbelian}
in~\S\ref{section:offenders}, after recalling Timmesfeld's replacement
theorem. In~\S\ref{section:centralSeries} we derive a result
(Lemma~\ref{lemma:Kn-2}) about offenders and central series
and apply this to the case of class at most four.
We treat the metabelian case
in~\S\ref{section:abelianDerived}, and show in \S\ref{section:MaxlClass}
that every group of maximal class has $p$-rank at most~$p$.
Finally, \S\ref{section:Rank-p}~is concerned with the case of a finite
$p$-group of $p$-rank at most~$p$ of Theorem~\ref{thm:main2}.

\section{Two notational conventions}
\label{section:notation}
\noindent
Throughout the paper, $G$ denotes a finite $p$-group, for an odd prime $p$.
We adopt the usual notation and conventions from the group theory literature
(see for example~\cite{Chermak:QuadraticAction}).
In addition, if $V$ is a \emph{right} $\f G$-module, we also view $V$ as an elementary
abelian $p$-group, written additively, and form the semidirect product
$\Gamma := G \ltimes V$, with $V$ the normal subgroup.
The group multiplication in~$\Gamma$ is
\[
(g,v)(h,w) = (gh, v * h + w) \, ,
\]
where we use~$*$ to denote the module action of $\f G$ on~$V$.

Our second convention regards the commutators. For
any $a,b\in G$, we set 
$$[a,b] := a^{-1} b^{-1} ab\quad\hbox{and}\quad[a_1,a_2,\ldots,a_{n+1}] = [[a_1,\ldots,a_n],a_{n+1}]~.$$
Inductively, we define $[a,b;1] = [a,b]$ and $[a,b;n] = [[a,b;n-1],b]$,
for all $n\geq2$.

\noindent
In particular, $g \in G$ identifies with $(g,0) \in \Gamma$ and
$v \in V$ with $(1,v) \in \Gamma$. Hence, the commutator
$[v,g]$ can be written as
\begin{align*}
[(1,v),(g,0)] & = (1,-v)(g^{-1},0)(1,v)(g,0) \\
& = (1,-v)(1,v*g) = (1,v*g-v) = (1,v*(g-1)) \, .
\end{align*}
That is,
\begin{equation}
[v,g] = v * (g-1) \, .
\label{eqn:commutator}
\end{equation}
In particular, $[v,g^p] = v * (g^p - 1) = v * (g-1)^p = [v,g;p]$, and so
\begin{equation}
[v,g^p] = [v,g;p] \, .
\label{eqn:pfold}
\end{equation}

\noindent
By identifying $G$ and $V$ with the corresponding subgroups of~$\Gamma$, one
obtains
\[
C_V(H) = \{v \in V \mid \forall h \in H \; v*h = v \} = V^H
\]
for a subgroup $H \leq G$, and
\[
C_G(V) = \{g \in G \mid \forall v \in V \; v*g = v \}  \trianglelefteq G \, .
\]
Note that the $\f G$-module $V$ is faithful if and only
if $C_G(V) = 1$. In this paper we will only be interested in faithful
modules.

\begin{defn}
Let $G$ be a $p$-group and $V$ a faithful $\f G$-module. A non identity
element $g \in G$ is called \emph{quadratic on~$V$\/}
if $[V,g,g] = 0$. If there is no confusion for~$V$, we simply say that
$g$ is \emph{quadratic}.
Since $[V,g^p] = [V,g;p]$, faithfulness implies that
quadratic elements must have order~$p$.
\end{defn}

\section{A lemma of Meierfrankenfeld and Stellmacher}
\label{section:MeierStell}

\begin{defn}[\cite{MeierfStellm:OtherPGV} 2.3]
Let $G$ be a finite group and $V$ a faithful $\f G$-module. For a subgroup
$H \leq G$ one sets
\[
j_H(V) := \frac{\abs{H}\abs{C_V(H)}}{\abs{V}} \in \mathbb{Q} \, .
\]
Note that $j_1(V)=1$.
\end{defn}

\begin{lemma}[Lemma 2.6 of~\cite{MeierfStellm:OtherPGV}]
\label{lemma:MeierfStellm}
Let $A$ be an abelian group and $V$ a faithful $\f A$-module. Let
$H,K$ be two subgroups of~$A$. Then
\[
j_{HK}(V) j_{H \cap K}(V) \geq j_H(V) j_K(V) \, ,
\]
with equality if and only if $C_V(H \cap K) = C_V(H) + C_V(K)$.
\end{lemma}

\begin{proof}
Since $\langle H,K\rangle=HK$, we have an equality $C_V(HK)=C_V(H)\cap
C_V(K)$. In addition, $C_V(H \cap K) \supseteq
C_V(H)+C_V(K)$. Therefore,
$$\begin{array}{lll}j_{HK}(V)j_{H\cap K}(V)&=\frac{\abs{HK}\abs{C_V(HK)}\abs{H\cap
    K}\abs{C_V(H\cap K)}}{\abs V^2}&\geq\\
&\geq\frac{\abs{HK}\abs{H\cap
    K}\abs{C_V(H)}\abs{C_V(K)}}{\abs V^2}&=j_H(V)j_K(V)~,\end{array}$$
because $\abs{HK} \abs{H\cap K} = \abs{H} \abs{K}$.
\end{proof}

\section{Offenders and abelian normal subgroups}
\label{section:offenders}

\begin{notn}
Let $G$ be a finite group and $p$ a fixed prime number.
We denote by $\E(G)$ the poset
of non-trivial elementary abelian $p$-subgroups of~$G$.
\end{notn}

\begin{defn}[\cite{RevisionVol2} 26.5]
Let $G$ be a finite group and $V$ a faithful $\f G$-module.
A subgroup $E \in \E(G)$ is an \emph{offender of ~$G$
on $V$} if $j_E(V) \geq 1$. If there is no confusion on $G$ and $V$, we
simply say that $E$ is an \emph{offender\/}. If $V$ has an offender,
then $V$ is called an \emph{$F$-module\/}. An offender $E$ is
\emph{quadratic on $V$\/} (or simply \emph{quadratic\/}) if $[V,E,E] = 0$.
Define
\[
\pp := \{E \leq G \mid \text{$E \in \E(G)$ and
$j_E(V) \geq j_F(V)$ $\forall\;1 \leq F \leq E$}\} \, .
\]
Consequently every $E \in \pp$ is an offender,
and every minimal\footnote{That is, a minimal element of the set of offenders,
partially ordered by inclusion.} offender lies in~$\pp$.
\end{defn}

\noindent
Note that $V$~is an $F$-module if and only
if $\pp$ is nonempty. The subgroups in $\pp$ are sometimes called
\emph{best offenders}.

\vspace{.3cm}
\noindent
We shall assume that the reader is familiar with Chermak's
treatment~\cite{Chermak:QuadraticAction}
of Timmesfeld's replacement theorem~\cite{Timmesfeld:Replacement}\@.

\begin{lemma}
\label{lemma:quadraticOffender}
Let $V$ be a faithful $\f G$-module and $E \in \pp$. Then there is a
quadratic offender $F \in \pp$
which satisfies $j_F(V)=j_E(V)$ and $F \leq E$.
\end{lemma}

\begin{proof}
This is Timmesfeld's replacement theorem
(\cite[Theorem 2]{Chermak:QuadraticAction})\@, applied with $F = C_E([V,E])$. The construction of~$F$ implies $[V,F,F]=0$.
\end{proof}

\begin{rk}
Note that Timmesfeld's replacement theorem also gives the decomposition
$\;C_V(F) = [V,E] + C_V(E)\;$.
\end{rk}

\begin{lemma}
\label{lemma:Enormal}
Let $G$ be a $p$-group and $V$ a faithful $\f G$-module such that there are
no quadratic elements in $\Omega_1(Z(G))$. Then $E \cap Z(G) = 1$
for every quadratic offender~$E$. In particular, a quadratic offender
does not contain any non-trivial normal subgroup of~$G$.
\end{lemma}

\begin{proof}
Let $E$ be a quadratic offender. We have $E \cap Z(G) = E \cap
\Omega_1(Z(G))$. Now, every non-trivial
element of~$E$ is quadratic, whereas no element of $\Omega_1(Z(G))$ is. Hence, $E \cap Z(G) = 1$.
The last statement follows from the fact that every non-trivial normal subgroup meets $Z(G)$.
\end{proof}

\noindent
We are now ready to show Theorem~\ref{thm:normalAbelian}.

\begin{proof}[Proof of Theorem~\ref{thm:normalAbelian}]
We first show the implication (\ref{enum:Enormal})${}\Rightarrow{}$(\ref{enum:Gdash}).
Let $E$ be an offender. Then $Z(G'E)$ is an abelian normal subgroup of~$G$, since $G'E \trianglelefteq G$.
By part~(\ref{enum:Enormal}), this means that $E \nleq Z(G'E)$.
Since $E$~is abelian, this can only happen if $[G',E] \neq 1$.

For part~(\ref{enum:Enormal}), suppose that~$A$ does contain an offender
$E$. Note that $E$ lies in the elementary abelian subgroup
$C := \Omega_1(A)$. Of course, $V$~is faithful as an $\f C$-module.
Set
\[
j_0 := \max \{ j_E(V) \mid E\in\E(C)\}
\]
Choose $E \in\E(C)$ with $j_E(V)=j_0$. Note that
every such~$E$ lies in~$\pp$. By Lemma~\ref{lemma:quadraticOffender}
we may assume that~$E$ is quadratic. Among the quadratic offenders
in~$A$ (and hence in $C$)
with $j_E(V)=j_0$, let us pick~$E$ of minimal order.
We claim that ~$E$ is a T.I. subgroup of $G$, that is $E\cap E^g=1$, for
all $g\in G-N_G(E)$. Note that since $C$ is normal in $G$, then all the
$G$-conjugates of $E$ lie in $C$. Observe also that $E$ lies in $\mathscr{M}(C,V)$, in the terminology of
\cite[Lemma~1]{Chermak:QuadraticAction}, since $C$~is elementary abelian
and we chose~$j_E(V)$ maximal. So every $G$-conjugate of~$E$ lies
in $\mathscr{M}(C,V)$.
By Lemma~1 of~\cite{Chermak:QuadraticAction},
the intersection of any family
of conjugates of~$E$ lies in $\mathscr{M}(C,V)$, and hence, the minimality assumption on the order of $E$ forces any such
intersection to be trivial, whenever it is a proper subgroup of $E$. So $E$ is a T.I. subgroup of $G$, as claimed.

Now $E \ntrianglelefteq G$ by Lemma~\ref{lemma:Enormal}\@. Hence $1$ is an
intersection of conjugates of~$E$. This implies that $j_0 = j_1(V)=1$. Let
$F$ be a $G$-conjugate of~$E$, with $F\ne E$, and hence $F \cap E = 1$,
since $E$~is a T.I. subgroup of $G$. Moreover, we have equalities $j_E(V)=j_F(V)=j_1(V)=1$,
and also $j_{EF} (V) \leq j_0 = 1$, as $EF \leq C$ is elementary abelian.
So from Lemma~\ref{lemma:MeierfStellm} we deduce
that $j_{EF}(V)=1$ and that $V = C_V(1) = C_V(E) + C_V(F)$.
So Lemma~\ref{lemma:EF} shows that $[V,E] \subseteq C_V(H)$, where
$H \leq A$ is the normal closure of $E$~in $G$.
The same argument shows that $[V,F] \subseteq C_V(H)$ for every
conjugate $F$~of $E$. We therefore deduce that
\begin{equation}
[V,H] \subseteq C_V(H) \, .
\label{eqn:VHandCVH}
\end{equation}
Now, $H\in\E(C)$, since $E\in\E(C)$ and $C$ is a normal
elementary abelian subgroup of $G$. So $H$ is itself a quadratic
offender, by Eqn.~\eqref{eqn:VHandCVH}\@. But $H$~is also normal
in~$G$, so $H$ contradicts Lemma~\ref{lemma:Enormal}\@.
\end{proof}

\begin{lemma}
\label{lemma:EF}
Let $A$ be an abelian subgroup of a $p$-group~$G$ such that $[A,B]=1$ for
every $G$-conjugate $B$~of $A$. Let $V$ be an $\f G$-module such
that $[V,A,A]=0$ and $V = C_V(A) + C_V(B)$ for every conjugate
$B \neq A$ of~$A$. Then $[V,A] \leq C_V(H)$,
where $H$ is the normal closure of $A$~in $G$.
\end{lemma}

\begin{proof}
It suffices to show that $[V,A] \leq C_V(B)$ for every $G$-conjugate
$B$~of $A$, since $C_V(H)$ is the intersection of all the $C_V(B)$.
For $B=A$, we have $[V,A,A]=0$, by assumption. For $B\neq A$, let $v \in V$,
$a \in A$ and $b \in B$. By hypothesis, there is a decomposition $v = u + w$
with $u \in C_V(A)$ and $w \in C_V(B)$. Thus,
\[
[v,a,b] = [u,a,b] + [w,a,b] = [w,a,b] = w*(a-1)(b-1)
= w*(b-1)(a-1) \, ,
\]
and so $[v,a,b] = [w,b,a]=0$, as was left to be shown.
\end{proof}

\section{Central series}
\label{section:centralSeries}
\noindent
Recall the following terminology. Given a finite
group $G$, the ascending central series is defined inductively
by $Z_0(G) = 1$ and $Z_{r+1}(G)$ is the normal subgroup of $G$
containing $Z_r(G)$ and such that $Z_{r+1}(G)/Z_r(G)=Z(G/Z_r(G))$, for
all $r\geq1$. The descending
central series is given by $K_1(G)=G$ and inductively $K_{r+1}(G) = [K_r(G),G]$, for
all $r\geq1$. The class $n$~of $G$ is the smallest number such that $Z_n(G)=G$.
This is also the smallest number such that $K_{n+1}(G)=1$. Note that
if the class is~$n$, then
\begin{equation}
K_{n+1-r}(G) \leq Z_r(G) \quad \text{($0 \leq r \leq n$).}
\label{eqn:centralSeries}
\end{equation}
Further details are given in \cite[III.1]{Huppert:I}, 
or also \cite[(8.7)]{Aschbacher:book}, where $K_r(G)$ is denoted by~$L_r(G)$.
From the Three Subgroups Lemma
(\cite[(8.7)]{Aschbacher:book}), we have
\begin{equation}
[K_r(G), K_s(G)] \leq K_{r+s}(G) \, .
\label{eqn:Kcomm}
\end{equation}
We shall make repeated use of the following lemma from~\cite{oliver}:
\begin{quote}
\noindent
\textbf{\cite{oliver}, Lemma 4.1.} \;
Suppose that $p$ is an odd prime, that $G\neq1$ is a finite $p$-group, and
that $V$~is a faithful $\f G$-module.
Suppose that $A,B \in G$ are such that $C := [A,B]$ is a non-trivial element
of $C_G(A,B)$. If $C$ is non-quadratic, then so are $A$~and $B$.
\end{quote}

\begin{lemma}
\label{lemma:Kn-2}
Let $G$~be a $p$-group of class~$n$\@, let $V$ be a faithful
$\f G$-module, and suppose that $E$ is a quadratic offender. 
If $2r \geq n$ and $K_{r+1}(G)$ contains no quadratic elements, then
\[
[K_r(G), E] = 1 \, .
\]
In particular if $n \geq 4$ and there are no quadratic elements in
$\Omega_1(Z(G))$, then
\[
[K_{n-2}(G), E] = 1 \, .
\]
\end{lemma}

\begin{rk}
The last part has no meaning for $n \leq 2$, and the example
discussed in \cite[\S5]{oliver} shows that it is false for $n=3$.
\end{rk}

\begin{proof}
We prove the first part by induction, starting with $r = n$ and working
downwards. We have $[K_n(G),E]\leq K_{n+1}(G)=1$, by Equation
(\ref{eqn:Kcomm}), and so the claim holds for $r=n$. Now, let $\frac
n2\leq r<n$ and
suppose that $[K_{r+1}(G), E] = 1$.
If $[K_r(G),E] \neq 1$, then there are $a \in K_r(G)$ and
$e \in E$ with $c := [a,e] \neq 1$. Since $c \in K_{r+1}(G)$, we have
$[c,a] \in [K_{r+1}(G),K_r(G)] \leq K_{2r+1}(G) \leq K_{n+1}(G)=1$.
Moreover,  $[c,e]=1$ since $c \in K_{r+1}(G)$ and the inductive hypothesis states
that $[K_{r+1}(G),E]=1$.
In other words, we have $1 \neq c = [a,e]$ with $e$ quadratic and
$[c,a]=[c,e]=1$. But then, \cite[Lemma 4.1]{oliver} says that $c$~is
quadratic, which contradicts
the assumption that no quadratic element lies in $K_{r+1}(G)$.

For the second part, if $n \geq 4$, we have $2r \geq n$ for $r =
n-2$. So, it is enough to show that $K_{n-1}(G)$ contains no quadratic
element. As recalled above, there is an inclusion $K_{n-1}(G) \leq
Z_2(G)$. By \cite[Lemma 4.1]{oliver}, since 
$\Omega_1(Z(G))$ does not contain any quadratic elements, there are no
quadratic elements in $Z_2(G)$ either.
\end{proof}

\begin{theorem}
\label{thm:class4}
Suppose that $G$ is a $p$-group and $V$ is a faithful $\f G$-module such that
$\Omega_1(Z(G))$ has no quadratic elements. If
$G$ has class at most four then $V$~cannot be an $F$-module.
\end{theorem}

\begin{proof}
Recall that $V$ is an $F$-module if and only if $\pp$ is not empty. 
By definition of $\pp$, every offender contains an element of~$\pp$, and
by Timmesfeld's replacement theorem (Lemma~\ref{lemma:quadraticOffender}),
every offender contains a quadratic offender which lies in~$\pp$.
It therefore suffices to prove that there are no quadratic offenders.
Let $E\in\E(G)$. From part~(\ref{enum:Gdash}) of Theorem~\ref{thm:normalAbelian}, we have
$[G',E] \neq 1$ if $E$ is a quadratic offender.

If $G$ has class three or less, then $[G',E] \leq Z(G)$. Since
$\Omega_1(Z(G))$ contains no quadratic elements, neither does the setwise commutator $[G',E]$. 
Thus, \cite[Lemma 4.1]{oliver} shows that $E$ is not a quadratic offender.

Now assume the class is four.
In this case we have
$G' = K_2(G) = K_{4-2}(G)$ and so $[G',E]=1$ by
Lemma~\ref{lemma:Kn-2}\@. So, $E$ is not a quadratic offender.
\end{proof}

\section{Metabelian groups}
\label{section:abelianDerived}
\noindent
Recall that a group is metabelian if and only if its derived subgroup is
abelian.
In this section we will need the following well-known result:

\begin{lemma}
\label{lemma:Huppert}
Suppose that $G$ is a finite group and that $A$ is an
abelian normal subgroup such that $G/A$ is cyclic, generated by the coset
$xA$ of $x \in G$. Then
$G' = \{[a,x] \mid a \in A\}$.
\end{lemma}

\begin{proof}
This is Lemma~4.6 of Isaacs' book~\cite{Isaacs}, or more precisely the
equality $\theta(A)=G'$ in the proof. It is also
Aufgabe~2\,a) on p.~259 of Huppert's book~\cite{Huppert:I}\@.
\end{proof}

\begin{theorem}
\label{thm:abelianDerived}
Suppose that $G$ is a $p$-group and $V$ is a faithful $\f G$-module such that
$\Omega_1(Z(G))$ has no quadratic elements. If the derived subgroup~$G'$
is abelian, then $V$~cannot be an $F$-module.
\end{theorem}

\begin{proof}
As noted in the proof of Theorem \ref{thm:class4}, if $V$~is an $F$-module then there is a quadratic offender~$E$.
From Part (\ref{enum:Gdash}) of Theorem~\ref{thm:normalAbelian}, we then
have $[G',E] \neq 1$. So, there is an $a \in E$ with $[G',a] \neq 1$.
The subgroup $K := G' \langle a \rangle$ of $G$ is a non-abelian normal subgroup of~$G$.
In particular, $K'\cap \Omega_1(Z(G))>1$. Let $c\in K'\cap\Omega_1(Z(G))$ be
a non-identity element. By Lemma~\ref{lemma:Huppert} there is an element $b \in G'$ with
$c = [b,a]$. So $c$ must be quadratic by \cite[Lemma~4.1]{oliver},
since $a$~is quadratic and $c$~is central.
But that cannot be, for $c$~lies in $\Omega_1(Z(G))$ and is therefore
non-quadratic by assumption.
\end{proof}

\begin{corollary}
\label{coroll:abelianCyclic}
Suppose that $G$ is a $p$-group and $V$ is a faithful $\f G$-module such that
$\Omega_1(Z(G))$ has no quadratic elements. Suppose that $G$ has
an abelian normal subgroup $A \trianglelefteq G$ such that $G/A$ is abelian too.
Then $V$~cannot be an $F$-module.
\end{corollary}

\begin{proof}
As $G/A$ is abelian, we have $G' \leq A$ and hence $G'$ is abelian.
\end{proof}

\section{Maximal class}
\label{section:MaxlClass}

\begin{lemma}
\label{lemma:maxlClassRankp}
Let $G$ be a finite $p$-group of maximal class. Then the $p$-rank of~$G$ is at
most~$p$. Moreover, if $p$~is odd, only the wreath product
$C_p \wr C_p$ has $p$-rank $p$.
\end{lemma}

\begin{rk}
This fact does not appear to be generally known. However, we believe that it is
remarked in passing by Berkovich~\cite{Berkovich:subgroupsEpimorphic}\@.

Let us also point out that the definition of $p$-groups of maximal
class may vary. Indeed, in~\cite{Huppert:I}, Huppert allows abelian groups of
order $p^2$ to be of maximal class, whereas in~\cite{LeedhamGreenMcKay:book}
Leedham-Green and McKay stipulate that groups of maximal class have
order at least~$p^4$. We follow Huppert's conventions.
\end{rk}

\begin{proof}
Let $G$ be a finite $p$-group of order $p^n$ and maximal class $n-1$.
First we consider the small cases with $n \leq p+1$.
An abelian group of order~$p^2$ has $p$-rank at most two.
A nonabelian group of order~$8$ can have $2$-rank at most two.
If $p$~is odd and $G$ is nonabelian with an elementary abelian subgroup of
rank~$p$, then its order must be at least $p^{p+1}$. 
So to finish off the cases with $n \leq p+1$, we just need to consider
the case where $p$~is odd, $n=p+1$ and $G$~contains an elementary
abelian subgroup $V$~of rank~$p$. As $G$~is not abelian,
the $p$-rank of $G$~is $p$. As $V$~has index $p$~in $G$,
it is normal
and the factor group $G/V$ has order $p$. 
Since $G$ has class~$p$, the group $G/V$
acts on $V$ as one $(p \times p)$-Jordan block with eigenvalue~$1$. 
Let $a\in G\setminus V$. Then, $a$ acts on $V$ with minimal polynomial
$(x-1)^p$ and $G=\langle V, a\rangle$. 
Note that $a^p \in V$ lies in the one-dimensional eigenspace, i.e. the
center $Z(G)$ of $G$.
So replacing $a$~by $ab$ for a suitably chosen element $b$~of the
set-theoretic difference $V \setminus [V,a]$, we may assume that
$a^p = 1$.
Hence the extension splits, and $G$ is isomorphic to the wreath product
$C_p\wr C_p$.

From now on we assume that $\abs{G} = p^n$ with $n \geq p+2$, and we
appeal to the following results of \cite[III]{Huppert:I}. 
By 14.16 Satz, we have that $G_1 = C_G(K_2(G)/K_4(G))$ is regular, and that
$\abs{G_1 \colon \mho_1(G_1)} = p^{p - 1}$, and so, by 10.7 Satz,
$\abs{\Omega_1(G)} = p^{p-1}$. 
Thus, if $G$ contains an elementary abelian subgroup $V$ of rank $p$,
then $V$ contains some $g \in G \setminus G_1$. It is hence enough to
show that no element of order $p$ in $G\setminus G_1$ has a
centralizer of $p$-rank greater than $p-1$.
By 14.6 Hauptsatz b), $G$ is not
exceptional, as $n \geq p+2$. Recall that by 14.5 Definition, exceptional groups arise
for $n\geq 5$, whence this result is relevant only for $n\geq p+2$ and
for $p$ odd. Now, we apply 14.13 Hilfsatz b), which states that if
$g \in G \setminus G_1$, then $\abs{C_G(g)} = p^2$. Consequently,
$G$ has at most rank $p-1$.
\end{proof}

\section{Rank \emph{p}}
\label{section:Rank-p}

\begin{theorem}
\label{theorem:rank-p}
Suppose that $p$~is odd, that $G$ is a $p$-group, and that $V$ is a faithful
$\f G$-module such that every non-identity element of $\Omega_1(Z(G))$ acts with minimal polynomial $X^p-1$.
If the $p$-rank of $G$ is at most~$p$, then $V$~cannot be an $F$-module.
\end{theorem}

Recall that the condition that every element of $\Omega_1(Z(G))$ acts
with minimal polynomial $X^p-1$ is condition (PS) of~\cite{oliver}.

\begin{proof}
Suppose that $V$ is an $F$-module. By Lemma~\ref{lemma:quadraticOffender},
$\pp$ contains a quadratic offender~$E$. Since $E \cap Z(G) = 1$ and $G$~has
rank~$p$, the rank of~$E$ can be at most $p-1$.
So by Lemma~\ref{lemma:preRankp} below, the rank of~$E$ is exactly~$p-1$.
Moreover, the normal closure $F$~of $E$~in $N_G(N_G(E))$ is elementary
abelian of rank~$p$, since it is strictly larger than~$E$.
Pick $h \in N_G(N_G(E)) \setminus N_G(E)$.
Then $E^h \neq E$, so $F = \langle E,E^h\rangle$. Moreover
$E\cap E^h$ must have size at least $p^{p-2}$ and is therefore nontrivial.
Pick $1 \neq c \in E^h \cap E$. By
Lemma~\ref{lemma:preRankp} Eqn.~\eqref{eqn:CVg}, this means that
$C_V(E)=C_V(c)=C_V(E^h)$ and therefore $C_V(F)=C_V(E)$. So given
that $j_E(V)=1$, the definition of~$j_F$ means that $j_F(V)=p$.
But this contradicts Timmesfeld's replacement theorem
(Lemma~\ref{lemma:quadraticOffender}),
as there is no quadratic offender~$H$ with $j_H(V)=p$.
\end{proof}

\begin{lemma}
\label{lemma:preRankp}
Suppose that $p$~is odd, that $G$ is a $p$-group, and that $V$ is a faithful
$\f G$-module such that every non-identity element of $\Omega_1(Z(G))$ acts with minimal
polynomial $X^p-1$.
\begin{enumerate}
\item
If $E \in \E(G)$ is an offender, then its rank is at least $p-1$.
\item
If $E \in \E(G)$ is a rank $p-1$ offender, then $E$~is a quadratic
offender and lies in~$\pp$. Moreover, $j_E(V)=1$; we have
\begin{equation}
\label{eqn:CVg}
C_V(g)=C_V(E) \quad \text{for every $1 \neq g \in E$;}
\end{equation}
and the normal closure
$F$~of $E$ in $N_G(N_G(E))$ is elementary abelian with $F > E$.
\end{enumerate}
\end{lemma}

\begin{proof}
Every offender contains an element of~$\pp$, which in turn contains
a quadratic offender. So minimal elements of the poset of offenders
are quadratic and lie in~$\pp$.
So it is enough to consider the case of a quadratic
offender~$E$.

For all $g\in E$ the subspace $C_V(g)$ of $V$ is $Z(G)$-invariant:
if $z \in Z(G)$ and $v \in C_V(g)$ then
\[
[v*z,g] = v*z*(g-1)=v*(g-1)*z=0 \, .
\]
Now choose a non-identity element $g\in E$, and let $i$ be the smallest
integer such that $[V,Z(G);i] \subseteq C_V(g)$.
We claim that $[V,Z(G);i] \subseteq C_V(g^h)$ for every $h \in G$.
To see this, note that $[V,Z(G);i]$ is an invariant subspace of~$V$,
and therefore $[V,Z(G);i] * h^{-1} = [V,Z(G);i]$. Hence
\[
[[V,Z(G);i],g^h] = [V,Z(G);i]*h^{-1} * (g-1) * h =  [V,Z(G);i] * (g-1) * h
= 0 \, .
\]
This means that $[V,Z(G);i] \subseteq C_V(H)$, where $H$ is the normal
closure of $\langle g \rangle$~in $G$. 
So $H\cap\Omega_1(Z(G))\neq1$, and for any $z\in H\cap\Omega_1(Z(G))$
we have that $[V,z;i+1]=0$.
The minimal polynomial assumption on $\Omega_1(Z(G))$ therefore implies that
$i+1 \geq p$: and so the definition of~$i$ means that
$\abs{V:C_V(g)}\geq p^{p-1}$.
As $C_V(E) \subseteq C_V(g)$, it follows that $\abs{V:C_V(E)} \geq p^{p-1}$.
But $\abs{E} = j_E(V) \abs{V:C_V(E)}$, and $j_E(V) \geq 1$ since
$E$~is an offender.
This proves the first part.

Second part:
If $E$ is a rank $p-1$ offender, then it is minimal in the poset of offenders,
therefore quadratic and a member of~$\pp$. The proof of the first part
shows Eqn.~\eqref{eqn:CVg}\@.
It also shows
that $j_E(V)=1$ and that $\abs{V:C_V(E)} = p^{p-1}$.

Set $N = N_G(E)$. Observe that $N_G(N) > N$, since $E \ntrianglelefteq G$
by Lemma~\ref{lemma:Enormal}\@. This shows that $F > E$.
Pick $g,h \in N_G(N)$ with $Ng \neq Nh$.
Then $E^g,E^h \leq N$; and as $N_G(E^g) = N_G(E)^g = N$, we see
that $E^g,E^h$ normalize each other. If they always centralize
each other, then $F$ is indeed elementary abelian, since it is generated
by all such~$E^g$. So suppose $a \in E^g$, $b \in E^h$ have
nontrivial commutator $c = [a,b]$. As $E^g,E^h$ normalize each other, we have
$c \in E^g \cap E^h$. Since $E^g,E^h$ are rank~$p-1$ offenders,
it follows by Eqn.~\eqref{eqn:CVg} that $C_V(a)=C_V(c)=C_V(b)$. So
$[a,b]=1$ by Lemma~\ref{lemma:commonCentralizer} below, a contradiction.
\end{proof}

\begin{lemma}
\label{lemma:commonCentralizer}
Let $G$ be a $p$-group and $V$ is a faithful $\f G$-module.
Suppose that $a,b \in G$ satisfy either of the following conditions:
\begin{enumerate}
\item
$a,b$ are quadratic, and $C_V(a) =C_V(b)$; or
\item
$[V,a] \subseteq C_V(b)$ and $[V,b] \subseteq C_V(a)$.
\end{enumerate}
Then $[a,b]=1$.
\end{lemma}

\begin{proof}
Recall that $[V,a]\subseteq C_V(a)$, for any quadratic element $a\in
G$. So $[V,a] \subseteq C_V(b)$ and $[V,b] \subseteq C_V(a)$
do hold if $a,b$ are quadratic and have the same centralizer in~$V$.
That is, the first case is a consequence of the second one.

In the second case we have
$[V,a,b]=[V,b,a]=0$, using additive notation in the $\f
G$-module $V$. Hence, a routine computation yields $[V,[a,b]]=0$. Since
$V$ is faithful, we deduce that $[a,b]=1$, as claimed.
\end{proof}

\end{document}